\documentclass[12pt]{article}
\usepackage{cite}
\usepackage{authblk}
\usepackage{lineno}
\usepackage{caption}
\usepackage{amsmath}
\usepackage{wrapfig}
   \usepackage{amssymb}
   \usepackage{graphicx}
   \usepackage{amscd}
   \usepackage{rotate}
   \usepackage{amsfonts}
   \usepackage{mathrsfs}

\modulolinenumbers[5]
\newtheorem{definition}{Definition}
\newtheorem{theorem}{Theorem}
\begin{document}

%\begin{frontmatter}

\title{On a stochastic gene expression with pre-mRNA, mRNA and protein contribution}

\author[1]{Ryszard Rudnicki}
\author[2]{Andrzej Tomski}
\affil[1]{Institute of Mathematics, Polish Academy of Sciences, Bankowa 14, 40-007, Katowice, Poland ryszard.rudnicki@us.edu.pl}

\affil[2]{Institute of Mathematics, Jagiellonian University, \L ojasiewicza 6, 30-348 Krak\'{o}w, Poland
andrzej.tomski@im.uj.edu.pl}
\maketitle
\noindent DOI: 10.1016/j.jtbi.2015.09.012

\medskip

\noindent $<2015>$. This manuscript version is made available under the 

\noindent CC-BY-NC-ND 4.0 license

\begin{abstract}
In this paper we develop a model of stochastic gene expression, which is  an extension
of the  model investigated in the paper [T. Lipniacki, P. Paszek, A. Marciniak-Czochra, A.R. Brasier, M. Kimmel, Transcriptional stochasticity in gene expression, J. Theor. Biol. $238\ (2006)\ 348-367$].
In our model, stochastic effects still originate from random fluctuations in gene activity status, but we precede mRNA production by the formation of pre-mRNA, which enriches classical transcription phase. We obtain a stochastically regulated system of ordinary differential equations (ODEs) describing evolution of pre-mRNA, mRNA and protein levels.
We perform mathematical analysis of a long-time behaviour of this stochastic process, identified as a
piece-wise deterministic Markov process (PDMP). We check exact results using numerical simulations for the distributions of all three types of particles. Moreover, we investigate the deterministic (adiabatic) limit state of the process, when depending on parameters it can exhibit two specific types of behavior: bistability and the existence of the limit cycle. The latter one is not present when only two kinds of gene expression products are considered.

\end{abstract}


\section{Introduction}

Gene expression and its regulation is a very complex process, which takes place in the cells of living organisms,
especially in eukaryotes \cite{reed}. It is widely known that this process depends on the behaviour of
crucial substances, called transcription factors (TFs) and chromatin architecture. Our investigation is based on the idea of \cite{lip}, where a simplified diagram of gene expression was presented. It was mentioned there that genes fluctuate randomly between their activity or inactivity status and transcripts are produced in bursts. Stochastic effects at the initial stage are very strong compared to both the matter production and degradation processes, so we consider the noise of Markov-type origin merely at the activation stage. These claims were verified and analysed through the years \cite{blake}, \cite{friedman}, \cite{kepler}, \cite{kim}, \cite{pedr}. The whole scheme describes expression of a single gene, assuming
it has $n$ copies, but further analyse was performed in the case of one copy only. After activation of the gene (which is initiated by binding to the promoter region some of TFs), mRNA transcription and protein translation phases follow. At first, mature mRNA is produced in the nucleus, then it is transported from the nucleus to the cytoplasm, where the second
phase takes place. As a result, new proteins are born.

\begin{figure}[h]
\includegraphics[width=12cm]{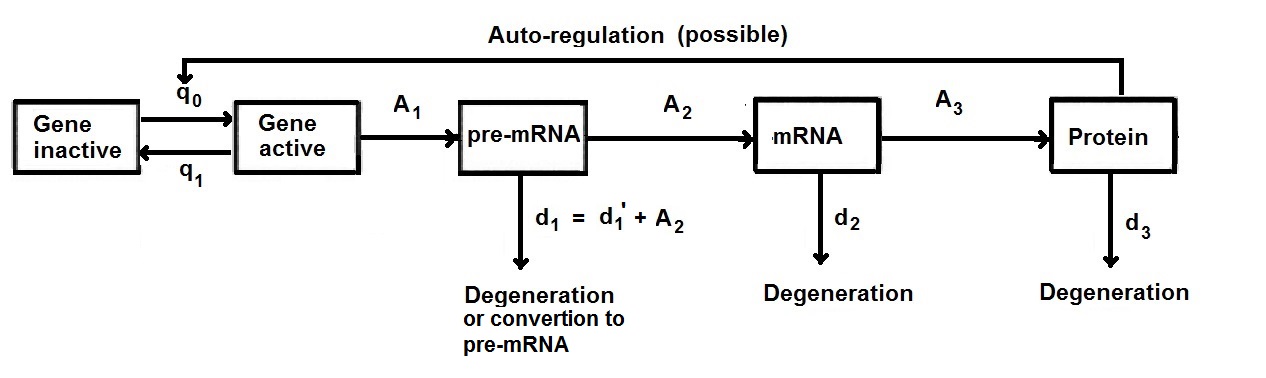}
\caption{An extended scheme of (auto-regulated) gene expression. }
\label{fig:1.1}
\end{figure}

     In the mentioned class of models, not only transcription and translation evolution were considered,
but also biological degradation of both types of the particles: mRNA and protein. All the processes were recognised as continuous, so the planar system of ordinary linear differential equations were used to represent the dynamics of fluctuations in the level of certain type particles. Moreover, first equation included stochastic ``switch'' component, being responsible for the control of gene activity status. This system has been identified in \cite{bobrow} as a
Piece-wise Deterministic Markov Process (PDMP), introduced by \cite{davi}. However, after reflection
on these results, an important question arises: to what extent does the two-stage model fits the current state of biological knowledge? Would adding another stage make description of the gene expression more precise? Finally, will the problem be much more complicated if we add the third stage? In the mentioned work of \cite{lip} there is
a remark that translated mRNA particle must get through some further processing
before a new, mature protein is formed. Beside that, plenty of thematic books; \cite{lodi}, \cite{wat} and publication sources; \cite{peng}, \cite{yap} claiming that at least one additional phase, called primary transcript (or pre-mRNA) processing should be taken into account. Actually, in eukaryotic genes, after the activation signal, the DNA
code is transformed into pre-mRNA form of transcript. Then, the non-coding sequences (introns) of transcript are cut off.
This action is combined with other modifications widely known as RNA processing. Only then we get a functional form of mRNA, which is transferred into the cytoplasm, where during the third phase, translation phase, mRNA is decoded into a protein.
In short, we consider three-phase model of gene expression with three main components, i.e. three
variables $x_1,x_2,x_3$ describing evolution of pre-mRNA, mRNA and protein levels. Firstly, we assume that pre-mRNA molecules
are produced at the rate $A_1 \gamma(t),$ where $A_1$ is a constant and we introduce a stochastic binary valued function $\gamma(t) \in \{0,1\} $ which marks, at time $t \ge 0$, if the gene is in active ($\gamma(t)=1$) or inactive ($\gamma(t)=0$) state. This function will be described in detail in Sec.~\ref{ss:23}. The mRNA production rate is equal to $A_2 x_1(t),$ where $A_2$ is a constant and $x_1(t)$ denotes the number of pre-mRNA molecules at time $t$. Similarly, the protein
translation takes the place at the rate $A_2 x_1(t),$ where $x_2(t)$ denotes the number of mRNA molecules at time $t.$
Moreover, all three types of particles undergo the degradation process. The total lost of pre-mRNA particles
is given by $d_1 x_1=d_1' x_1+A_2 x_1,$ where the constant $d_1'$ is the degradation rate of pre-mRNA particles and another constant
$A_2$ is the rate of converting pre-mRNA into mRNA particles. It means that $d_1=d_1' + A_2 $ should be treated as the
total degradation rate of pre-mRNA particles. This concept takes into consideration that pre-mRNA is converted to mRNA
\cite{lip2}, in contrast to mRNA which serves as a template for mRNA synthesis, but is not degraded during the
synthesis. Thus, in other cases we use standard description, i.e. the constants $d_2$ and $d_3$ denote, respectively, mRNA and protein degradation rates. This expansion of the previous, simplified diagram of gene
expression depicted in \cite{lip}, is now presented in Fig.~\ref{fig:1.1}. We note that the switching between active and inactive state of the gene depends on
the so-called jump rates (activation/inactivation rates).

Introduction of the third variable to the model means, that its geometry moves unavoidably into $\mathbb{R}^3$ space. Although,
in the last few years some PDMP-based biological models were presented, they focus on the applications in planar systems: \cite{bobrow}, \cite{lip}, \cite{przem}, \cite{tom}. What is worth mentioning, in the paper of \cite{bakhtin} it was investigated that three-dimensional Lorenz system with stochastic switching ``admits a robust strange attractor", but here we concentrate on a situation, when
jump rates are not necessarily constant and also on the convergence in time of the distribution of the process to the equilibrium distribution. In the proof of the main theorem, we use some results concerning asymptotic stability of Markov semigroups. The main idea is that we need to check two conditions: irreducilibity of the semigroup and the existence of some estimation of the semigroup from below which can be checked by using H\"ormander condition for considered process. An alternative (with no reference to Markov semigroups) approach concerning long-time behavior of PDMP is based on regularity and convergence results from the papers \cite{ben} and \cite{sben2}, where similar conditions for the convergence of PDMP were developed. Recently, \cite{mal} asked the question, is it possible to describe long-time qualitative properties of the process for spatial dynamics? Here we do such an analysis with the aim to include the role of primary transcript in the processes basic for eukaryotic genes.

      This paper is organised as follows. First, in Section~\ref{s:2} we present idea of the model. Then we discuss the deterministic (adiabatic) limit state, which can exhibit two specific types of behavior: bistability and the existence of the limit cycle. Later, we present mathematical description of the process, including deterministic part and the
      stochastic component. Having done that, we introduce Markov semigroups and we recall how
they can be generated by PDMP to describe time evolution of the densities of the process. In the first part of Section~\ref{s:3} we formulate the main theorem of this paper, which says that the Markov semigroup related to the model is asymptotically stable. This means that there exists a stationary density and independently on the initial distribution the density of the process converges to the stationary density as time goes to infinity. We find a set, an ``attractor", on which this three-dimensional distribution is concentrated. The second part of Section~\ref{s:3}  is devoted to stochastic simulations of the process. We show time-dependent and mutual dependent behaviour of levels of pre-mRNA, mRNA and protein. We also approximate the above-mentioned limit stationary density for all types of the particles. In the further part we also refer to the deterministic state behavior. In Section ~\ref{s:4} we sum up the results of our paper and give some conclusion remarks.
\section{The model}
\label{s:2}
\subsection{Construction}
\label{ss:21}
Let $x_1,x_2,x_3$ denote three non-negative variables, which describe time-evol\-ving levels of pre-mRNA, mRNA and protein,
respectively. In accordance with the current surveys, we consider that the activation and inactivation rate functions
$q_0(x_1,x_2,x_3)$ and $q_1(x_1,x_2,x_3)$ can be constant \cite{bakhtin}, \cite{peccoud} or can depend on the number of the particles of one type,
usually the proteins \cite{crudu}, \cite{lip}. Briefly speaking, the gene is activated with the rate $q_0(x_1,x_2,x_3)$ and inactivated with the rate $q_1(x_1,x_2,x_3)$.
The minimal mathematical assumptions about $q_0$ and $q_1$ in the case of two variables are discussed by \cite{bobrow}.
In line with the approach of \cite{lip}, we study evolution of the following system of ODEs with a stochastic component:
\begin{equation}
\label{e:1}
\begin{cases} 0 \xrightarrow {q_0(x_1,x_2,x_3)} 1, \  0 \xleftarrow {q_1(x_1,x_2,x_3)} 1 \\
\dfrac{dx_1}{dt} =A_1\gamma(t)-d_1 x_1 \\ \dfrac{dx_2}{dt} = A_2 x_1-d_2 x_2 \\ \dfrac{dx_3}{dt} = A_3 x_2- d_3x _3 ,
\end{cases}
\end{equation}
where $A_i$ and $d_i$ are positive constants.

{\it Remark $1.$} We pay attention
to the fact that a standard three-dimensional Goodwin model of an oscillatory gene regulation loop \cite{goo} was developed in a
similar manner to ours and it can be interpreted even in the same way \cite{wang1}. However, instead of the presence of the stochastic process $\gamma(t)$,
Goodwin model contains a non-linear term describing the production rate of mRNA. The source of nonlinearity is the dependence of this rate from the protein level.
In our work we take this fact into account by making the intensity functions $q_0$ and $q_1$ possibly dependent from the level of
any type of particles, especially the proteins.

 If $q_0$ and $q_1$ are constant, we calculate the expected levels of pre-mRNA, mRNA and protein in the molecular population:
\begin{align*}
\mathbb{E}(x_1)&=\frac{A_1 q_0}{d_1(q_0+q_1)},\\
\mathbb{E}(x_2)&=\frac{A_1 A_2 q_0}{d_1 d_2(q_0+q_1)},\\
\mathbb{E}(x_3)&=\frac{A_1 A_2 A_3 q_0}{d_1 d_2 d_3(q_0+q_1)},
\end{align*}
despite the fact that these levels oscillate in time (see Sec.~\ref{ss:32}  for details). Using standard rescaling techniques known from investigation
of the planar model in \cite{bobrow}, we obtain the system:
\begin{equation}
\label{e:2}
\begin{cases} 0 \xrightarrow {q_0(x_1,x_2,x_3)} 1, \  0 \xleftarrow {q_1(x_1,x_2,x_3)} 1 \\ \dfrac{dx_1}{dt} =\gamma(t)-x_1 \\ \dfrac{dx_2}{dt} = a(x_1-x_2) \\ \dfrac{dx_3}{dt} = b(x_2-x_3) ,\end{cases}
\end{equation}
$a,b>0$, in addition $a\ne b,$ $a\ne 1$ and $b\ne 1$ . We investigate this system in the next sections of the paper. Our results remain true
also if $a=b$ or $a=1$ or $b=1$ (see Remark $3.$)

{\it Remark $2$.} We can consider a more complicated process containing larger number of intermediate steps which lead to
equations like in the system \ref{e:2} and the approach presented below will not change. However, we analyse the system \ref{e:2} with three equations
for the brevity of notation. Larger number of intermediate steps introduce time delay, which in the case of
negative feedback makes the system oscillatory. We have observed such oscillatory behaviour even in the three-dimensional case but for
very special rate functions $q_0$ and $q_1$ (see Fig. ~\ref{fig:3.24}).

\subsection{The adiabatic limit}
\label{s:21}

We shall consider particularly interesting behavior of our model, when both of the jump rates $q_0$ and $q_1$ tend to infinity, unlike their ratio. In this case, we can replace the stochastic process $\gamma(t)$ by its
expected value $\Gamma:= \mathbb{E} \gamma = \frac{q_0}{q_0+q_1}$ \cite{bobrow2,lip} to obtain a state called
{\it deterministic} or {\it adiabatic} limit. Hence, the system $\ref{e:2}$ transforms to deterministic system of
three ODEs:
\begin{equation}
\label{e:222}
\begin{cases}  \dfrac{dx_1}{dt} =\Gamma(x_3)-x_1 \\ \dfrac{dx_2}{dt} = a(x_1-x_2) \\ \dfrac{dx_3}{dt} = b(x_2-x_3). \end{cases}
\end{equation}
Depending on the values of the parameters $a,b$ we investigate some specific types of behavior \cite{jar}.
Firstly,  we consider the case of the positive autoregulation, i.e. when  $\Gamma$  is an increasing function of $x_3$. Assume that the equation $\Gamma(c)=c$ has three roots $c_1<c_2<c_3$
in the interval $(0,1)$ and $\Gamma'(c_1)<1$, $\Gamma'(c_2)>1$, $\Gamma'(c_3)<1$. Then the system (\ref{e:222}) has three stationary points
$\mathbf x_i=(c_i,c_i,c_i)$, $i=1,2,3$.
The linearization  of (\ref{e:222}) at  $\mathbf x_i$ leads to the following characteristic polynomial
\[
P(\lambda) =(1+\lambda)(a+\lambda)(b+\lambda)-\Gamma'(c_i)ab.
 \]
 Since $\Gamma'(c_i)>1$ we have
 $P(0)<0$, and $P(\lambda)>0$ for sufficiently large, which means that
 the polynomial $P$ has a positive root and, therefore, the point
 $\mathbf x_2$ is unstable. Now we check that stationary points are
$\mathbf x_1$ and $\mathbf x_3$ are asymptotically stable, i.e. all roots of $P$
have negative real parts. Indeed, if $a=b=1$ then $P$ has all roots with negative real parts:
\[
\lambda_1=\Gamma'(c_i)^{1/3}-1,  \quad \lambda_2=-1+(-\tfrac12+\tfrac{\sqrt{3}\,i}{2})\Gamma'(c_i)^{1/3},
\quad\lambda_3=-1+(-\tfrac12-\tfrac{\sqrt{3}\,i}{2})\Gamma'(c_i)^{1/3}.
\]
If we find  some coefficients $a,b$ such that
$P$ has a root with a nonnegative real part, then we find some coefficients $a,b$
such that  $P$ has a root with  a zero real part, i.e., $\lambda=\alpha i$,
$\alpha\in \mathbb R$. But then $\alpha=0$ or
\[
(a+b+ab)=\alpha^2 \textrm{ and }(1+a+b)\alpha^2=(1-\Gamma'(c_i))ab.
 \]
Both cases are impossible if $\Gamma'(c_i)<1$, hence we obtain a bistable state (see \cite{jar}). In such a case one can expect
that the stationary density of the stochastic process will be bimodal provided that $q_0$ and $q_1$ are finite, but sufficiently large.
We check this by performing extensive numerical simulations presented in Sec. \ref{ss:32}
In the case of  the negative  autoregulation, i.e. when $\Gamma$  is a decreasing creasing function of $x_3$,
we have only one stationary point $\mathbf x=(c,c,c)$, where $c$ is the unique solution of the equation $\Gamma(c)=c$.
Observe that in the case $a=b=1$  the polynomial $P$ has one negative real root, and two complex roots with positive real parts if $\Gamma'(c)<-8$. This suggest that in this case the limit cycle can appear.
We check its existence by simulating the system \ref{e:222} in Sec. \ref{ss:32}. This case is especially interesting, since from the standard Bendixson$-$Dulac theorem it follows that such limit cycle oscillations are not observed in the two-dimensional system studied before.
Again, one can expect that for sufficiently large $q_0$ and $q_1$ the stationary density for the process will be distributed close to the limit cycle trajectory.

\subsection{Two deterministic systems}
\label{ss:22}

For a fixed state of the gene, which determines the value of $\gamma(t)\equiv i,\ i \in \{0,1\},$ the process is
purely deterministic and we get the system of the first order differential equations
\begin{equation}
\label{e:3}
\begin{cases}  \dfrac{dx_1}{dt} =i
-x_1 \\ \dfrac{dx_2}{dt} = a(x_1-x_2) \\ \dfrac{dx_3}{dt} = b(x_2-x_3),
 \end{cases}
\end{equation}
with the initial condition $\mathbf {x}_0=(x_1^0, x_2^0, x_3^0) \in \mathbb{R}^3_+ $ with $a,b>0,\ a\ne b, a\ne 1$ and  $b\ne 1.$
The solution
$\pi_i^t(\mathbf {x_0})$ of this system is
\begin{equation}
\pi_i^t(\mathbf {x}_0)=i\mathbf {1} +  \mathrm{exp\,}(Mt)(\mathbf {x_0} - i\mathbf {1}),
\end{equation}
where $\mathbf {1}=[1,1,1]$ and
\[
M=\left[ \begin{array}{ccc}
- 1 & 0 & 0 \\
a & -a & 0 \\
0 & b & -b \\
\end{array} \right].
\]
Moreover, with a similarity to the two-dimensional case \cite{bobrow}, we have:
\begin{equation}
\label{e:sim}
\pi_t^1\mathbf (\mathbf {x_0})=\mathbf 1-\pi_t^0\mathbf (\mathbf 1)+ \pi_t^0\mathbf (\mathbf {x_0}).
\end{equation}

\begin{figure}[ht]
\includegraphics[height=9cm,width=15cm]{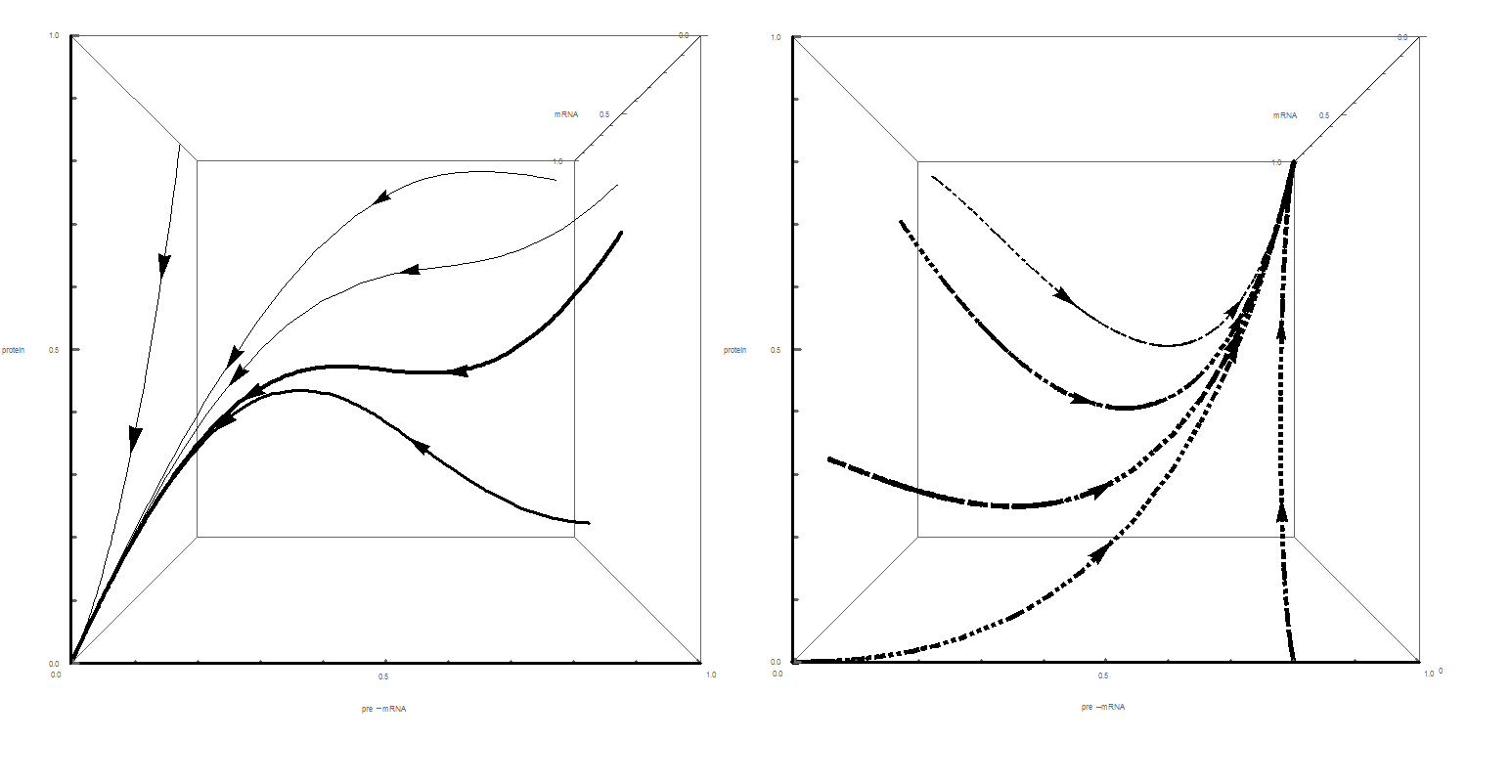}
\caption{A sample solutions of Eq. (\ref{e:3}) for $a=2,b=3,i=0$ (left)
and $a=2,b=3,i=1$ (right). }
\label{fig:2.2}
\end{figure}

In Fig.~\ref{fig:2.2} phase portraits of the system (\ref{e:3}) for both values of $i \in \{0,1\}$ are shown.
Each time, there exists one stationary solution: for $i=0;$ a point $(0,0,0)$ is asymptotically stable steady state,
as is a point $(1,1,1)$ for $i = 1.$ Looking at the right-hand sides of system \ref{e:3}, we state that if $x_1 > 1,$ then
(no matter what the value of $i$ is), $x_1$ decreases, moreover $x_2$ as well $x_3$ follow $x_1$. On the other hand, if $x_1 < 1,$ then it stays in the interval $[0, 1]$ forever, oscillating between $0$ and $1,$ the same happens with $x_2$ and $x_3.$ Hence, we can reduce the phase space of the process to a cube $X=[0,1]^3.$

\subsection{PDMP: a definition }
\label{ss:23}
We will briefly mention an idea behind PDMP introduced in ~\cite{davi}. We consider $q_0(x_1,x_2,x_3)$ and
$q_1(x_1,x_2,x_3) $ as two continuous and non-negative functions on $\mathbb{R}^3 $ such that:
\[
 q_0\left(0,0,0\right) \neq 0
\textrm{ and }q_1\left(1,1,1\right) \neq 0.
\]
Let $i_0 \in\{0,1\},\ T_0=0,\ \mathbf {x}_0 \in \mathbb{R}_+^3 $ and we define
a (random) function
$\gamma\colon [0,\infty) \to \{0,1\}$ satisfying
$\gamma(0)=i_0$
and
 \begin{eqnarray}
\gamma(t) :=
 \left\{
 \begin{array}{ll}
       i, &\textrm{if} \ T_{n-1}\leqslant t<T_{n},\\
       1-i, &\textrm{if} \ t=T_{n},
\end{array}
\right.
\end{eqnarray}
where for $n \geqslant 1,$ $T_n$ is a positive random variable
satisfying:
\begin{equation}
\begin{aligned}
F_{\mathbf {x}_s}(t)&=
\operatorname{Prob}(T_{n} - T_{n-1} \leqslant t|\ \gamma(T_{n-1})=i) \\
&= 1 - \operatorname{\mathrm{exp\,}} \left( - \int_{0}^{t} q_{i} (\pi_{i}(s,\mathbf {x}_{s}))ds \right),
\end{aligned}
\end{equation}

\begin{eqnarray}
\mathbf {x}_{s}:=
 \left\{
 \begin{array}{ll}

\pi^{s-T_{n-1}}_i(\mathbf {x}_{n-1}), &\textrm{if} \ T_{n-1}\leqslant s<T_{n},\\
     \mathbf {x}_{n}, &\textrm{if} \ s=T_{n},
\end{array}
\right.
\end{eqnarray}
where
\begin{equation}
\mathbf {x}_{n}=\pi^{T_n-T_{n-1}}_i(\mathbf {x}_{n-1}).
\end{equation}
In consequence, replacing a constant value $i \in \{0,1\}$ in the system $(\ref{e:3})$ by a stochastic
process $\gamma(t)$:
\begin{equation}
\label{e:4}
\begin{cases}  \dfrac{dx_1}{dt} = \gamma(t)
-x_1 \\ \dfrac{dx_2}{dt} = a(x_1-x_2) \\ \dfrac{dx_3}{dt} = b(x_2-x_3), \end{cases}
\end{equation}
gives a definition of a Markov process $\zeta(t)$ called a \textit{piece-wise deterministic Markov
process}, described by the quartet:
\begin{equation}
\zeta(t) := (x_1(t),x_2(t),x_3(t), \gamma(t)) =(\mathbf {x}(t),\gamma(t)).
\end{equation}
The state space of this process is $ \mathbb{X}=X \times \{0,1\}.$ The remaining characteristics
are the jump rates $q_i$ and the jump distribution $\mathbb{J}((\mathbf {x},i),\cdot)$ being the Dirac measure
$\delta_{(\mathbf {x},1-i)}$ such that
\begin{equation}
\mathbb{J}((x,i), \mathbb{X})=1.
\end{equation}
A random variable $T_n$ is called a \textit{time of the n-th jump of the process.} In ~\cite{bobrow} it was shown
that in such a case $\Delta_k=T_k-T_{k-1}>0,$ where $k \geqslant 1,\ \Delta_k < \infty$
and
\begin{equation}
\lim_{k \rightarrow \infty} T_k = \infty,
\end{equation} which means that the process is well-defined for all times $t \geq 0$.

\subsection{Markov semigroups and their link with PDMP}
\label{ss:24}
Now we will recall some definitions about Markov semigroups. We use them to describe the evolution of
distributions of the process given by the system (\ref{e:2}).
Detailed information about some connections
between semigroup theory and stochastic processes can be found in \cite{la} or \cite{rud2000}.
Let $(\mathbb{X},\Sigma,m)$ be a $\sigma-$finite measure space and let $D \subset L^1=L^1(\mathbb{X},\Sigma,m)$
be the set of the densities, i.e.
\[
 D =\{f \in L^1: f \geqslant 0,\  ||f||=1 \}.
\]
\begin{definition}
A linear $D$ preserving mapping $P: L^1 \rightarrow L^1$ is called a \textit{Markov} (or \textit{stochastic}) \textit{operator}.
\end{definition}
\begin{definition}
A family $\{P(t)\}_{t \geqslant 0} $ of Markov operators, which satisfies the following conditions:
\begin{itemize}
\item $P(0)=$Id  (identity condition),
\item $P(t+s)=P(t)P(s) $ for $s,\ t \geqslant 0$ (semigroup condition),
\item for each $f \in L^1$ the function $t \rightarrow P(t) f $ is continuous with respect to the $L^1$ norm  (strong continuity),
\end{itemize}
is called a \textit{Markov semigroup}.
\end{definition}
\begin{definition}
 A Markov semigroup $\{P(t)\}_{t \geqslant 0} $ is \textit{partially integral} if there exist
 $t_0 >0$ and a measurable function $k: \mathbb{X} \times \mathbb{X} \rightarrow \mathbb{R}^+,$ such that
 for every $f \in D:$
\begin{equation}
 \int_{\mathbb{X}} \int_{\mathbb{X}} k(p,q) m(dp)m(dq) > 0
\end{equation}
and
\begin{equation}
P(t_0)f(p) \geqslant \int_{\mathbb{X}} k(p,q) f(q) m(dq).
\end{equation}
\end{definition}
\begin{definition}
A Markov semigroup $\{P(t)\}_{t \geqslant 0}$ is \textit{asymptotically stable} if
\begin{itemize} \item there exists an \textit{invariant density} for $\{P(t)\}_{t \geqslant 0}$, i.e. $f^* \in D$ such that $P(t) f^* =f^*$ for all $t>0,$
\item for every density $f \in D:$
\end{itemize}
\begin{equation}
\lim_{t \rightarrow \infty} ||P(t)f-f^*||=0.
\end{equation}
\end{definition}
Below we define a property, which is in some sense ``opposite" to asymptotic stability, introduced in \cite{kom}.
\begin{definition}
A Markov semigroup is \textit{sweeping} (or \textit{zero-type}) with respect to a set $A \in \Sigma$ if for every $f \in D:$
\begin{equation}
\lim_{t \rightarrow \infty} \int_A P(t)f(x)m(dx)=0.
\end{equation}
\end{definition}

A precise instruction on how to construct Markov semigroup for PDMP is given by \cite{bobrow}.
Using the analogy with the two-dimensional model, we write Fokker-Planck system of equations for the partial densities
$f_0,f_1$ of the process
\begin{equation}
\label{e:5}
\begin{cases}
\dfrac{\partial f_0}{\partial t} + \dfrac{\partial }{\partial x_1}(-x_1 f_0)
+ a \dfrac{\partial }{\partial x_2}\left((x_1-x_2) f_0\right) + b \dfrac{\partial }{\partial x_3}\left((x_2-x_3) f_0\right)=
q_1 f_1 - q_0 f_0 \\\dfrac{\partial f_1}{\partial t} + \dfrac{\partial }{\partial x_1}\left((1-x_1) f_1\right)
+ a \dfrac{\partial }{\partial x_2}\left((x_1-x_2) f_1\right) + b \dfrac{\partial }{\partial x_3}\left((x_2-x_3) f_1\right)=
q_0 f_0 - q_1 f_1 ,
\end{cases}
\end{equation}
where $f_0,f_1$ are the functions defined on $[0,\infty) \times [0,1]^3$ such that for any Borel set
$\mathfrak{B} \subset \mathbb{R}^+ \times \mathbb{R}^+ \times \mathbb{R}^+$
\begin{equation}
\operatorname{Prob}(x(t) \in \mathfrak{B},\, \gamma(t)=i) =\iiint\limits_
{\mathfrak{B}} f_i(t,x_1,x_2,x_3)\,dx_1\,dx_2\,dx_3,\, i=0,1.
\end{equation}
For the reason of the presence of three spatial variables and a wide range of possible jump rates,
system (\ref{e:5}) is difficult to be solved analytically.
However, we will use Markov semigroup $\{P(t)_{t \geq 0}\}$ generated
by this process to prove that it has stationary density, which is an equilibrium with respect to time evolution of the distributions.

\section{Results}
\label{s:3}
\subsection{Asymptotic stability}
\label{ss:31}
\begin{figure}
\includegraphics[scale=0.4]{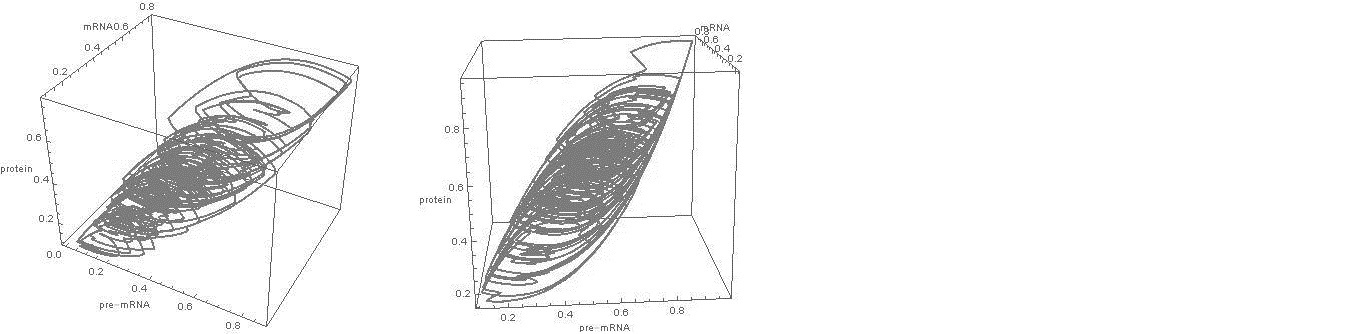}
 \centering
\caption{The sketch of $A$ obtained by numerically portrayed trajectories of the process in the three-dimensional space. The initial condition is
$x(0)=\frac{1}{2}, y(0)=\frac{1}{2}, z(0)=\frac{1}{2}, \gamma(0)=0$ and $a=2, b=10$ for constant (left) and protein-mediated (right) jump rates.}
\label{fig:3.11}
\end{figure}

In this section we present main result of this paper. We consider two particular solutions of the system (\ref{e:2}).
The first, $\phi(t)=(\phi_1(t),\phi_2(t),\phi_3(t))$
is the solution of (\ref{e:3}) with $i=0$ and the initial condition $(\phi_1(0),\phi_2(0),\phi_3(0))=(1,1,1)$.
The second $\psi(t)=(\psi_1(t),\psi_2(t),\psi_3(t))$
is the solution of
(\ref{e:2}) with $i=1$ and the initial condition $(\psi_1(0),\psi_2(0),\psi_3(0))=(0,0,0)$.
We conclude that $\phi$ and $\psi$ are two solutions of the system (\ref{e:2}) with $i=0$ and $i=1,$ respectively,
which join the asymptotically stable points $(0,0,0)$ and $(1,1,1)$
We construct the set $A$ in the following way.
Let $A_0$ be the surface made of all solutions of the system (\ref{e:2})
with $i=1,$ which start from any point lying on $\phi$.
This is also the case with $A_1,$ being the surface made of all the solutions of the system (\ref{e:2})
with $i=0$, which start from any point lying on $\psi$.
\begin{figure}[ht]
\begin{center}\includegraphics[height=9cm,width=11cm]{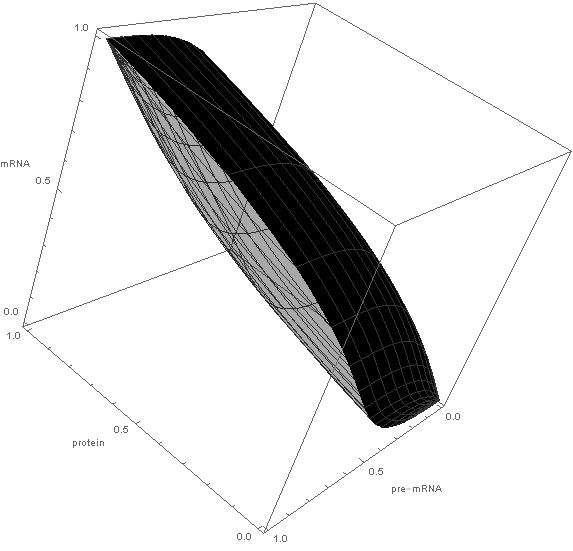}
\end{center}
\caption{The boundaries of $A,\ A_0 $ - filled and $A_1$ - transparent for $a=2$ and $b=10.$ }
\label{fig:3.1}
\end{figure}
We derive algebraic formulas describing $A_0$ as well as $A_1$ in Appendix A. Having done that, we define $A$ as a
subset of $[0,1]^3,$ bounded by $A_0$ and $A_1.$ In Fig.~\ref{fig:3.1} we show geometric visualisation of $A$. For comparison,
in Fig.~\ref{fig:3.11}, we portray the sketch of this set,  obtained by numerical simulations of the trajectories of the process (see Sec. 3.2).

 Now we can formulate the main result of the paper.
 \medskip
 \medskip

\begin{theorem}
Let $q_i(\mathbf {x})>0,\ i=0,1$ for $\mathbf {x} \in  [0,1]^3.$ Then, the Markov semigroup $\{P(t)\}_{t \geqslant 0}$ is
asymptotically stable and the support of the invariant density is the set $\mathbb{A}= A \times \{0,1\}$, where $A$ is expressed in the basis $(\ref{eigvec})$ by
\begin{equation}
A=\{(x-y+z,\ x^a-y^a+z^a,\ x^b-y^b+z^b): 1\geqslant x\geqslant y \geqslant z \geqslant 0 \}.
\end{equation}
\end{theorem}

A general idea beyond the strict proof of this theorem is provided by \cite{bobrow}. However, the proof in the two-dimensional model is simpler, because all the properties of the attractor are easy to deduce using geometrical arguments. For example, in this case the proof that $A$ is an invariant set for the process follows immediately from the M\"uller theorem \cite{walter} and the communication between states inside $A$ follows from the Darboux property.
Since the geometric arguments in the three-dimensional model are not obvious, we need to use precise formula to define the set $A$ and to prove its properties. Namely, we follow \cite{rud2} and prove that $\mathbb{A}$ is a set such that
\begin{itemize}
\item $\mathbb{A}$ is invariant for the process, i.e. if $(x_1(0),x_2(0),x_3(0),\gamma(0)) \in \mathbb{A},$
then $(x_1(t),x_2(t),x_3(t),\gamma(t)) \in \mathbb{A}$ for any $t>0,$
\item trajectories $(x_1(t),x_2(t),x_3(t),\gamma(t))$ of the process starting from any arbitrary point
from $[0,1]^3 \times \{0,1\}$ converge to $\mathbb{A}$ when time goes to infinity,
\item there is no smaller set satisfying these two conditions above.
\end{itemize}
Detailed mathematical proofs of these claims are long and provided in Appendix B. In Fig.~\ref{fig:4.1}
we show two-dimensional projections of A onto the 2D plane, looking exactly the same as the set proposed by \cite{lip}.
\begin{figure}[h]
\includegraphics[height=8cm,width=13cm]{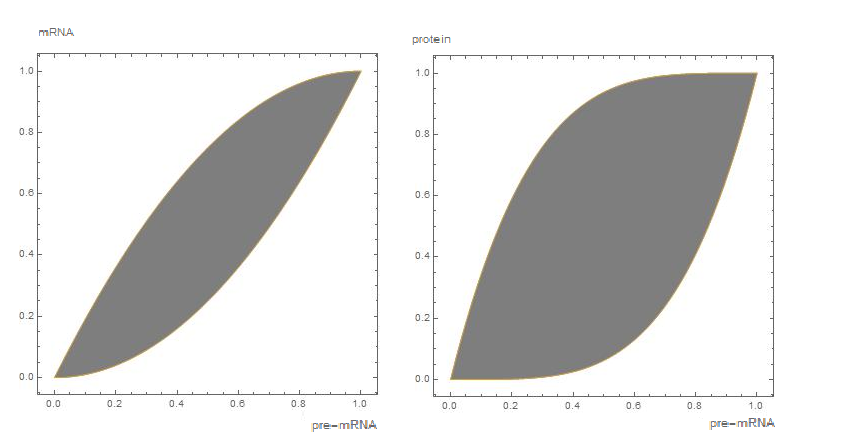}
\caption{
Projections of $A$ onto the plane, $a=2, b=3.$ Left: pre-mRNA and mRNA. Right: pre-mRNA and protein.
 }
\label{fig:4.1}
\end{figure}

\subsection{Stochastic simulations}
\label{ss:32}

Although it is difficult to solve Fokker-Planck equations $(\ref{e:5})$ analytically,
here we discuss stochastic simulations of the trajectories and distributions
of the system $(\ref{e:4})$, made to check the accuracy of our statements.
Such an approach, based on the \cite{gil} algorithm, was used by \cite{przem} to visualize
the evolution of the trajectories in the model of self-renewal cells differentiation.
For our model a similar code in Wolfram Mathematica environment was generated and run.

\begin{figure}
\includegraphics[scale=0.4]{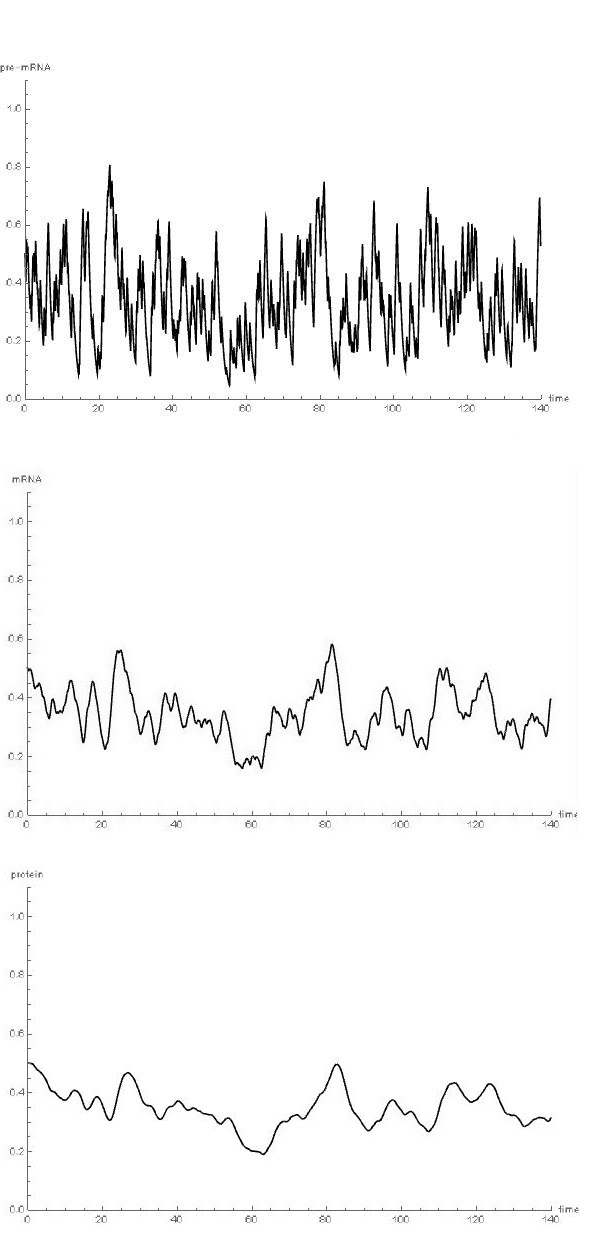}
\caption{Trajectories of the stochastic process $(\ref{e:4})$. The initial condition is
$x(0)=\frac{1}{2}, y(0)=\frac{1}{2}, z(0)=\frac{1}{2}, \gamma(0)=0$ and $a=\frac{1}{2}, b=\frac{1}{3},q_0=3, q_1=6 $ are set to show the level of pre-mRNA (top),
mRNA (center) and protein (bottom). }
\label{fig:3.2}
\end{figure}

\begin{figure}
\includegraphics[scale=0.4]{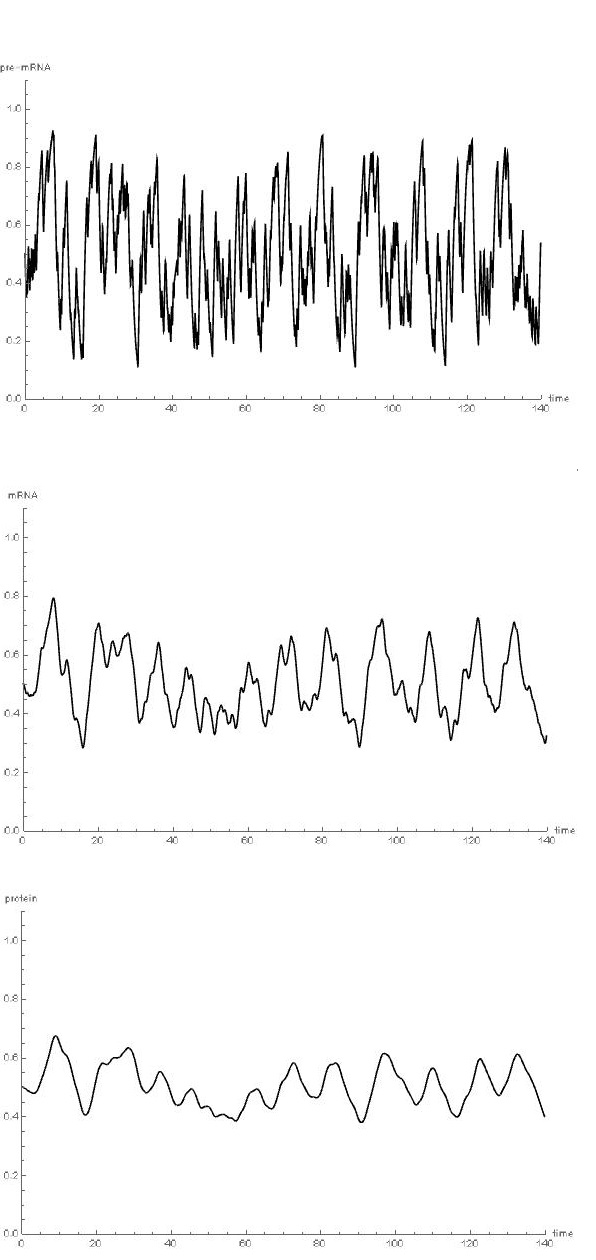}
\caption{Trajectories of the stochastic process $(\ref{e:4})$. The initial condition is
$x(0)=\frac{1}{2}, y(0)=\frac{1}{2}, z(0)=\frac{1}{2}, \gamma(0)=0$ and $a=\frac{1}{2}, b=\frac{1}{3},q_0=3, q_1=6\, z(t) $ are set to show the level of pre-mRNA (top),
mRNA (center) and protein (bottom). }
\label{fig:3.3}

\end{figure}

\begin{figure}
\includegraphics[scale=0.4]{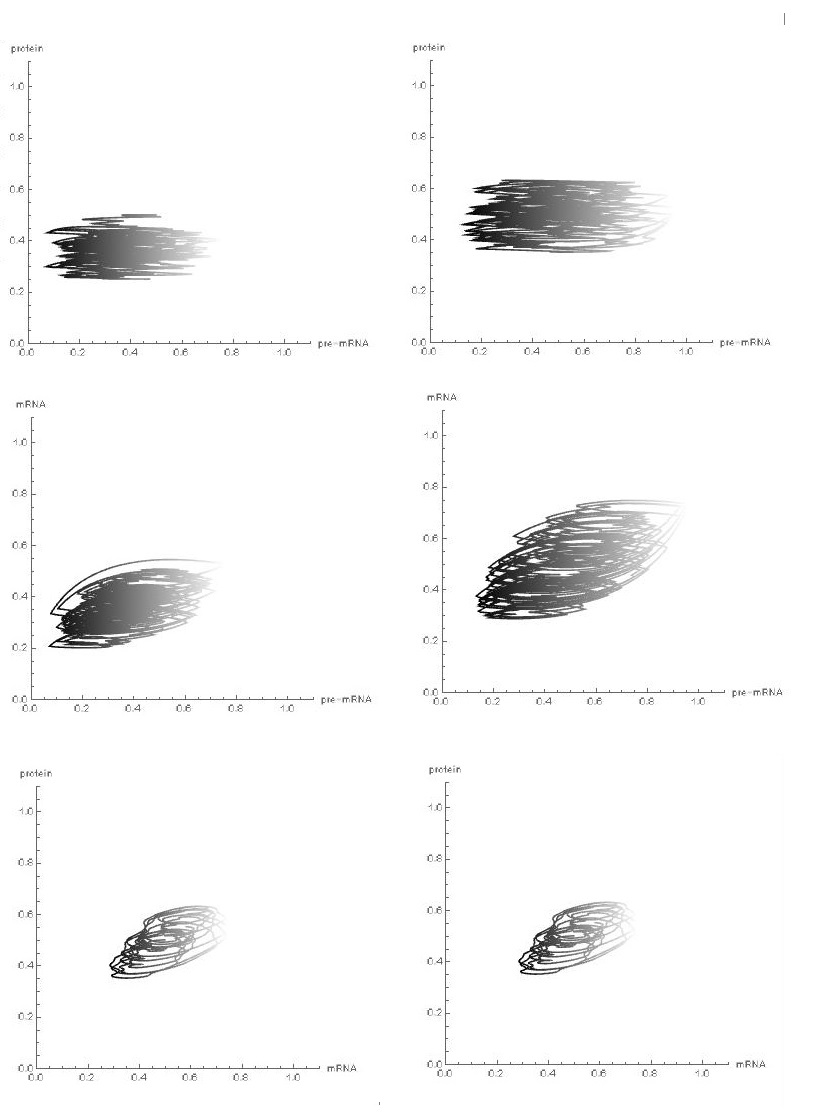}
\caption{Two-dimensional trajectories of the stochastic process $(\ref{e:4})$. The initial condition is
$x(0)=\frac{1}{2}, y(0)=\frac{1}{2}, z(0)=\frac{1}{2}, \gamma(0)=0$ and $a=\frac{1}{2}, b=\frac{1}{3}, q_0=3$  are set to show the dependence between pairs of the variables, for constant ($q_1=6$, left) and protein-mediated ($q_1=6 x_3,$ right) jump rates. }
\label{fig:3.23}
\end{figure}

We have compared the trajectories of the system $(\ref{e:4})$ for selected values of the parameters $a,b$
and jump rates $q_0(x_1,x_2,x_3),q_1(x_1,x_2,x_3)$ up to the final time moment $T=150$
or at least $1300$ jumps were performed. Afterwards, we depicted time evolution of the levels
of all three kinds of gene expression products for $a=\frac{1}{2}, b=\frac{1}{3}, q_0=3$
and $q_1=6 $ with the initial condition $x(0)=y(0)=z(0)=\frac{1}{2},$ separately in Fig.~\ref{fig:3.2}
and pairwise in Fig. ~\ref{fig:3.23} (left).
In Fig.~\ref{fig:3.3} and in Fig. ~\ref{fig:3.23} (right), where we use the same parameters as in the previous case,
we assume that the inactivation rate depends on the protein level, i.e. $q_1(x_1,x_2,x_3)=q_1 x_3.$

We notice that in both cases, as it was expected, fluctuation in pre-mRNA level is much stronger
than it is in mRNA or protein level. However, when the jump rates are constant, then
all three levels seem to vary in a more limited range than in the case when the jump rates are protein-mediated:
the values of standard deviation for consecutive phases are $0.151, 0.0887, 0.0609$ and $0.182, 0.099, 0.0617,$ respectively.
Moreover, we empirically calculated correlation level between each two phases. While pre-mRNA and mRNA levels
(with the values of the coefficient $0.555$ for constant and $0.573$ for protein-mediated jump rates), as well as mRNA
and protein ($0.685$ and $0.668$) levels were significantly correlated, pre-mRNA and protein levels were poorly related to each other
($0.07$ and $0.058$).

\begin{figure}
\includegraphics[scale=0.50]{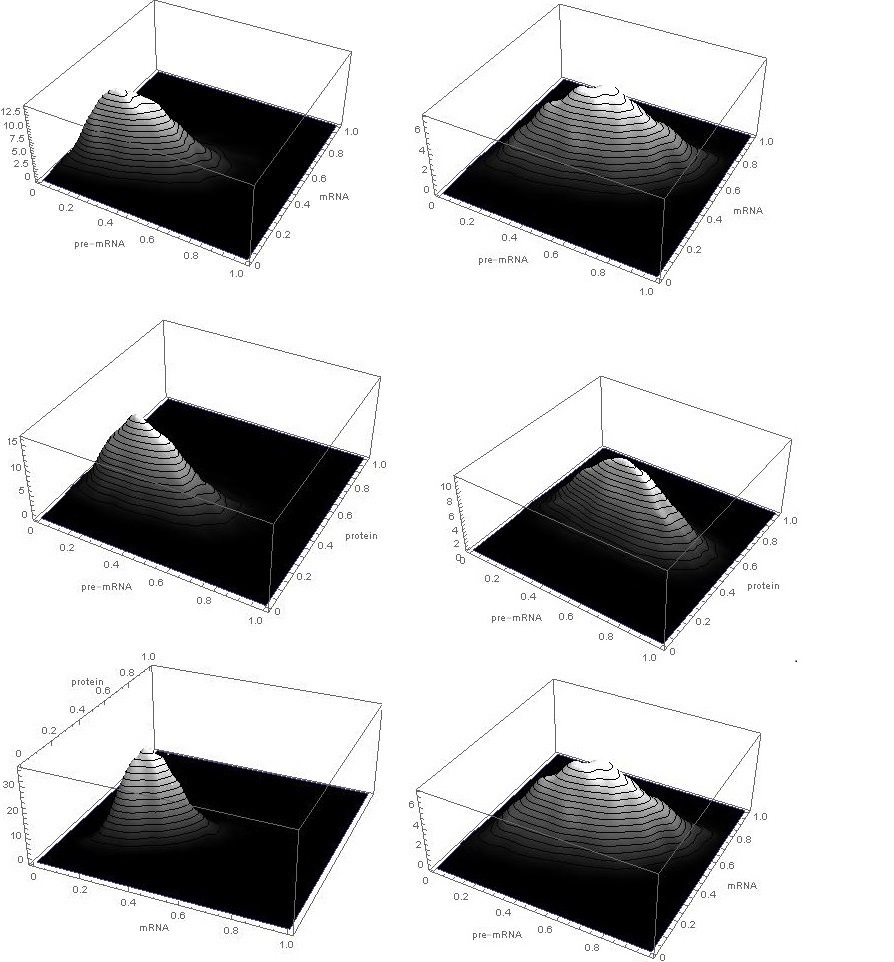}
\caption{Marginal distributions $\rho(t,x_k,x_j)=f_0(t,x_k,x_j)+f_1(t,x_k,x_j)$ calculated for time $t=15$ with simulations
of the system $(\ref{e:4})$ repeated $5000$ times for constant (left) and protein-mediated (right) jump rates. }
\label{fig:3.4}
\end{figure}

\begin{figure}
\includegraphics[scale=0.50]{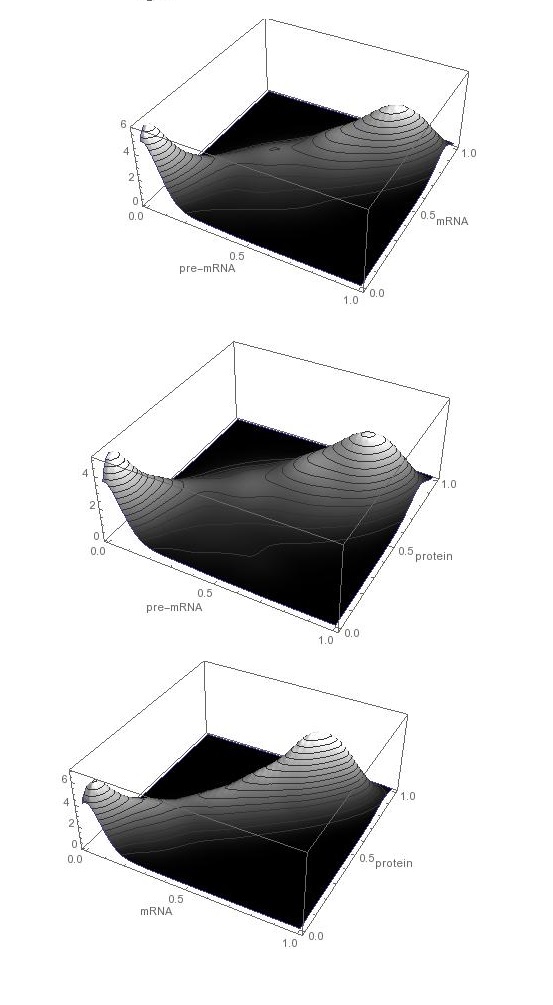}
\caption{Marginal bimodal distributions $\rho(t,x_k,x_j)=f_0(t,x_k,x_j)+f_1(t,x_k,x_j)$ calculated for $t=15$ with simulations
of the system $(\ref{e:4})$ with $a=b=1$ and jump rate functions $q_0(x_3)=10(0.01+{x_3}^2), q_1=10\cdot0,2.$ }
\label{fig:3.44}
\end{figure}

\begin{figure}
\includegraphics[scale=0.50]{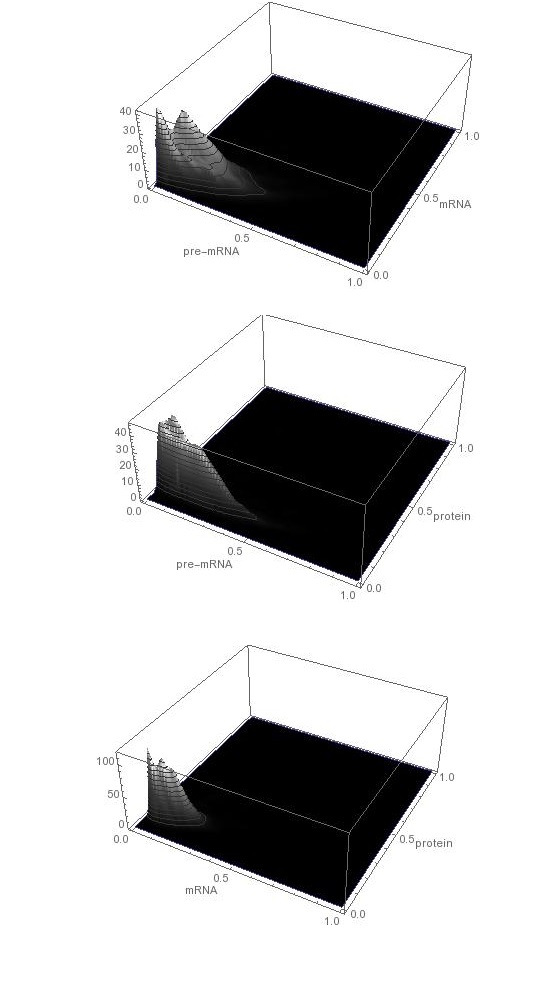}
\caption{Marginal distributions $\rho(t,x_k,x_j)=f_0(t,x_k,x_j)+f_1(t,x_k,x_j)$ calculated for $t=15$ with simulations
of the system $(\ref{e:4})$ with $a=b=1$ and jump rate functions$q_0=n$, $q_1(x_3)=9\cdot 10^{10}\,nx_{3}^{10}$  with $n=10.$ }
\label{fig:3.45}
\end{figure}

\begin{figure}

\includegraphics[scale=0.6]{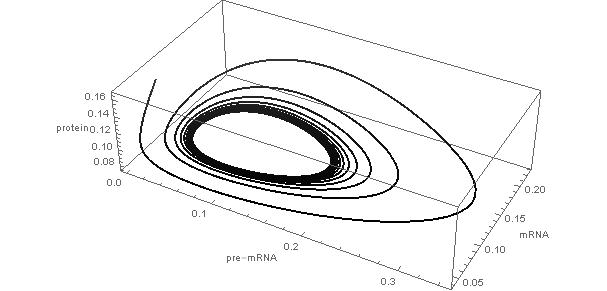}
\caption{The limit cycle for $a=b=1$ and jump rate functions $q_0^n\equiv n$ and $q_1^n(x_3)=9\cdot{10}^{10}\ n \ x_3^{10}$ in the adiabatic limit.}
\label{fig:3.24}
\end{figure}

Leaving all parameters unchanged, we analysed the distributions obtained by the simulations of system $(\ref{e:4})$
with constant and linearly dependent inactivation rate function $q_1(x_1,x_2,x_3)$,
respectively; see Fig.~\ref{fig:3.4}.  To follow the behaviour of a gene,
we pictured two-phase marginal distributions, i.e. $\rho(t,x_1, x_2)$, $\rho(t,x_2, x_3)$, $\rho(t,x_1,x_3)$,
where $\rho(t,x_k, x_j)=f_0(t,x_k,x_j)+f_1(t,x_k,x_j).$ The graphs were made by simulating
the system $(\ref{e:4})$ up to $T_F=15$, repeated $5000$ times and hence it is expected
(see also \cite{lip}) that the points $(x(T_F),y(T_F),z(T_F),i)$ obtained this way approximate stationary distributions $f_i(x_1,x_2,x_3)$.
Fluctuation strength, described by the jump rates, decides about the broadness of $\rho$: stronger fluctuations
give the broader distribution. However, on the left hand side of Fig.~\ref{fig:3.4}  the jump inactivation rate is twice as large as the jump activation rate,
so the inactive state becomes dominant. As a result, the distribution much more significantly points into the zero direction.

In Sec. \ref{s:21} we discussed and justified two specific types of behavior of the process in the deterministic (adiabatic) limit. In the first case
it is expected that when $q_0$ and $q_1$ are sufficiently large and $\mathbb{E} \gamma$ is an increasing function of the
protein level $x_3,$ we should obtain that the stationary density is bimodal, see Fig. ~\ref{fig:3.44}. On the other hand, when
both of the jump rates are still large, but this time $\mathbb{E} \gamma$ is a decreasing function of the
protein level $x_3,$ we observe the existence of the limit cycle in the deterministic limit. However, after simulating the process for
a sufficiently long time, we should obtain its density distributed close to this limit cycle. A comparison of the distribution of the stochastic process (Fig. ~\ref{fig:3.45}) and its
deterministic approximation (Fig. ~\ref{fig:3.24}) is provided. The latter case pays particular attention, because it is a kind of
behavior which is not present in the two-dimensional system from \cite{lip}.

\section{Conclusion}
\label{s:4}
We have studied a model of stochastic gene expression with the contribution of three main phases. Our investigation is based on the
two-dimensional model introduced by \cite{lip}, including: activation of the gene, mRNA transcription and
protein translation. Activity of the gene is regulated stochastically, namely by a Piece-wise Deterministic Markov Process,
\cite{bobrow}. In eukaryotes, where the transcript is produced in bursts, we can neglect other sources of stochasticity.
However, many reports: \cite{peng}, \cite{lodi}, \cite{wat}, \cite{yap} suggest that at least one additional phase,
i.e. pre-mRNA level regulation should be considered as well. This moves the state space of deterministic part of the process into $\mathbb{R}^3.$
We have analysed long-time behaviour of densities of the process. Using Markov semigroup techniques, we have shown that its distribution converges to equilibrium,
i.e. there exists a stationary density, such that independently on
the initial distribution its evolution is being stabilised with this density, when time goes to infinity.
Moreover, we have found a set, an ``attractor", which is a support for the equilibrium.
Using statistical approach, we visualised the trajectories of the process and approximated stationary distributions
in the cases of constant and protein-mediated jump rates. Moreover, we discussed two specific types of behavior
of the process: bistability and the existence of limit cycle trajectory. To summarize, we obtained qualitative and statistical results for the long-time evolution of three-phase kinetics of
the eukaryotic gene. We notice that the main result is in agreement with the one from the two-dimensional model.
This suggests that a sequence of gene transformations described by equations
including stochastic activation and only production and degradation processes,
does not have a significant influence on stabilizing long-time distribution of the product levels. However, our analysis
has also shown that it is not entirely true that the three-dimensional system has necessarily analogous dynamics to
the two-dimensional one studied before: the appearance of the limit cycle in our system, is not possible in the
planar model. Thus, larger number of intermediate steps in the case of negative feedback can make
the system oscillatory.

Nonetheless, another question is how to formulate more general approach,
when there is a need to analyse phases which cannot be described by linear ODEs and how their limit distribution and
the behavior of the trajectories will change.

\section*{Acknowledgements}

We thank P.R. Pa\'zdziorek for discussion and support. We are also grateful to the reviewers for their invaluable suggestions in improving the paper.
This paper was partially supported by the State Committee for Scientific Research (Poland) Grant No. 2014/13/B/ST1/00224 (RR).
%Here are two sample references:

\section*{Appendix A. Simplifying the system of ODEs and derivation of the formula for the attractor }
\label{sec:7}
Let us consider the system (\ref{e:2}) with a given value of $i,\ i=0,1$. Equivalently we can rewrite it as $\mathbf {x}'=M\mathbf {x}+c$, where
\[
M=\left[
\begin{array}{ccc}
- 1 & 0 & 0 \\
a & -a & 0 \\
0 & b & -b \\
\end{array}
\right]
\textrm{ and }
c=\left[
\begin{array}{c}
 1 \\
0 \\
0\\
\end{array} \right].
\]
We notice that $M$ has three distinct eigenvalues: $-1,-a,-b$ and
we can choose the eigenvectors $\mathbf v_1,\mathbf  v_2,\mathbf  v_3$ in such a way that the vector $\mathbf 1=[1,1,1]$ has also the coordinates $[1,1,1]$ in
the basis of the eigenvectors.
Precisely, the eigenvectors $\mathbf v_1,\mathbf v_2,\mathbf v_3$ are given by the formulas:
\begin{equation}
\label{eigvec}
\begin{aligned}
\mathbf v_1&=\Big[1,\frac{a}{a-1},\frac{ab}{(a-1)(b-1)}\Big],\\
\mathbf v_2&=\Big[0,\frac{-1}{a-1},\frac{b}{(a-1)(a-b)}\Big],\\
\mathbf v_3&=\Big[0,0,1-\frac{ab}{(a-1)(b-1)}-\frac{b}{(a-1)(a-b)}\Big].
 \end{aligned}
 \end{equation}
We transform $(\ref{e:3})$ by rewriting it in the new basis and we obtain new formulas for the system $(\ref{e:3})$:
\begin{equation}
\begin{cases}
 \frac{dx_1}{dt} = - x_1 + i \\
 \frac{dx_2}{dt} = - a x_2 +a i \\
 \frac{dx_3}{dt} = - b x_ 3 +b i
\end{cases}
\end{equation}
for $i=0,1$.
Let $\mathbf {x}=[x_1,x_2,x_3]^T$ be a column vector and let $\pi_t^i(\mathbf {x})$ denote the solution of $(\ref{e:2})$ at time $t$
with the initial condition $\mathbf {x}.$ We get
\begin{equation}
\label{app:01}
\pi_t^0(\mathbf {x})=(e^{-t}x_1,e^{-at}x_2, e^{-bt}x_3)
\end{equation}
and
\begin{equation}
\label{app:02}
\pi_t^1(\mathbf {x})= \mathbf{1} + \pi_t^0(\mathbf {x}) - \pi_t^0(\mathbf{1}),
\end{equation}
where $\textbf{1}$ denotes now a column vector $[1,1,1]^T.$
By alternate compositing of these functions we have the formulas
\begin{equation}
\pi_{t_2}^1\pi_{t_1}^0(\mathbf {x})= \mathbf{1} + \pi_{t_1+t_2}^0(\mathbf {x}) - \pi_{t_2}^0 \mathbf{1}
\end{equation}
 as well
\begin{equation}
\pi_{t_2}^0\pi_{t_1}^1(\mathbf {x})= \pi_{t_2}^0 \mathbf{1} + \pi_{t_1+t_2}^0(\mathbf {x}) - \pi_{t_1+t_2}^0 \mathbf{1}
\end{equation}
for any $t_1,t_2>0.$
Now, substituting $e^{-t_2}:=\alpha,$ $e^{-(t_1+t_2)}:=\beta$ we get

\begin{equation}
\label{e:100}
\pi_{t_2}^0\pi_{t_1}^1(\mathbf {x})= (\alpha+\beta (x_1-1),\alpha^a+\beta^a (x_2-1),\alpha^b+\beta^b (x_3-1)),
\end{equation}
and
\begin{equation}
\label{e:200}
\pi_{t_2}^1\pi_{t_1}^0(\mathbf {x})= (1-\alpha+\beta x_1,1-\alpha^a+\beta^a x_2,1-\alpha^b+\beta^b x_3)
\end{equation}
where $1 \geq \alpha \geq \beta \geq 0$.
Taking as the initial points $\mathbf {x}=(0,0,0)$ in the formula $(\ref{e:100})$ and $\mathbf {x}=(1,1,1)$  in the formula $(\ref{e:200}),$ we get parametric equations for the surfaces $A_0$ and $A_1$ which are the boundaries of $A$ (see Sec.~\ref{s:3}):
\begin{align*}
A_0&=\{(\alpha-\beta,\alpha^{a}-\beta^{a},\alpha^{b}-\beta^{b}):\ 1 \geq \alpha \geq \beta \geq 0 \},\\
A_1&=\{(1-\alpha+\beta,1-\alpha^{a}+\beta^{a},1-\alpha^{b}+\beta^{b}):\ 1 \geq \alpha \geq \beta \geq 0 \}.
\end{align*}

\noindent Similarly:

\begin{align*}
\pi_{t_3}^0\pi_{t_2}^1\pi_{t_1}^0(\mathbf {x})&= (\alpha-\beta +\gamma x_1,\alpha^a-\beta^a + \gamma x_2,
\alpha^b-\beta^b + \gamma^b x_3),\\
\pi_{t_3}^1\pi_{t_2}^0\pi_{t_1}^1(\mathbf {x})&= (1-\alpha+\beta +\gamma(x_1-1) ,1-\alpha^a+\beta^a +\gamma^a (x_2-1),1-\alpha^b+\beta^b + \gamma^b (x_3-1)),
\end{align*}

\noindent where $1 \geq \alpha \geq \beta \geq \gamma \geq 0$.

Now, let $V=\{(x,y,z)\colon 0<z<y<x<1\}$ and we define the function $\mathbf f\colon V\to \mathbb R^3$ by the formula
\begin{equation}
\mathbf f(x,y,z)=(x-y+z,x^a-y^a+z^a,x^b-y^b+z^b).
\end{equation}
It is easy to check that $\mathbf f$ is a local diffeomorphism. In particular,
\begin{equation}
\label{e:300}
\mathbf f(V)= \{\mathbf x= \pi_{t_3}^0\pi_{t_2}^1\pi_{t_1}^0\mathbf{1}:\,\, t_1>0, \, t_2>0,\, t_3>0\}
\end{equation}
is an open set and $\mathbf f(V)$ is the interior of $A$. Hence
\begin{equation}
A=\{(x-y+z,\ x^a-y^a+z^a,\ x^b-y^b+z^b): 1\geqslant x\geqslant y \geqslant z \geqslant 0 \}
\end{equation}

\noindent and $A_0$ i $A_1$ are indeed the boundaries of $A.$ Moreover $A_1$ is the symmetrical image of $A_0$ (and vice versa) with respect to a point $\left(\frac{1}{2},\frac{1}{2}, \frac{1}{2} \right).$ From this property we have the equivalent definition of $A,$ i.e.:
\begin{equation}
\label{e:500}
A=\{(1-(x-y+z),\ 1-(x^a-y^a+z^a),\ 1-(x^b-y^b+z^b): 1\geqslant x\geqslant y \geqslant z \geqslant 0 \}.
\end{equation}

\section*{Appendix B. The proof of asymptotic stability}
\label{sec:8}
Now we will prove the main result of this paper. We use the following theorem.
\begin{theorem}
\label{th:2}
Let $\mathbb{X}$ be a compact metric space and $\Sigma$ be the Borel $\sigma-$algebra.
If a Markov semigroup $\{P(t)\}_{t \geqslant 0} $ satisfies two conditions:
\begin{itemize}
\item[\rm{(a)}] for every density $f$ we have $ \ \int_{0}^{\infty} P(t) f \,dt > 0  $ a.e.,
\item[\rm{(b)}] for every $q_0 \in \mathbb{X} $ there exist $\kappa >0$, $t>0$ and a measurable function
$\eta \geqslant 0$ such that $ \int \eta(p)\, m(dp) > 0 $ and
\[
P(t) f(p) \geqslant \eta(p) \int_{B(q_0,\kappa)} f(q) m (dq),
\]
for $p \in \mathbb{X}$, where $B(q_0,\kappa)$ is the open ball with center $q_0$ and radius $\kappa,$
\end{itemize}
then the semigroup $\{P(t)\}_{t \geqslant 0} $ is asymptotically stable.
\end{theorem}
 Theorem \ref{th:2} was formulated as Corollary $1$
 in \cite{bobrow} and follows from two earlier results:

\begin{theorem}\cite{pich-rud1}
\label{th:as}
If $\{P(t)\}_{t\ge 0}$ is a partially integral Markov
semigroup and has a unique
invariant density $f_*>0$,   then the semigroup $\{P(t)\}_{t\ge 0}$ is asymptotically stable.
\end{theorem}
\begin{theorem}\cite{rud3}
\label{th:sw}
If $\{P(t)\}_{t\ge 0}$ is a  Markov semigroup on a metric space, satisfies (a) and (b),
and has no invariant density,
then it is sweeping from compact sets.
\end{theorem}
From Theorem \ref{th:sw} it follows that if the space $\mathbb{X}$ is compact
and the semigroup satisfies (a) and (b), then it has a unique and positive invariant density.
Now Theorem~\ref{th:2}
follows immediately
from Theorem~\ref{th:as}.

Hence, the idea of the proof is as follows.
First, we check that all trajectories enter the set
$\mathbb{A}$ and this set is invariant.
This allows us to reduce the proof of asymptotic
stability only to the semigroup $\{P(t)\}_{t \geq 0}$ restricted to the set $\mathbb{A}$.
Then we check that conditions (a) and (b)  are satisfied on $\mathbb{A}$.
Since $\mathbb{A}$ is compact the semigroup $\{P(t)\}_{t \geq 0}$  is asymptotically stable.

Firstly, we introduce some necessary definitions.
\begin{definition}
Let $V(M)$ be the set of real smooth vector fields on the
manifold $M$ on $ \mathbb{R}^d$ and let $C^{\infty}(M)$ denote the set of a real-valued smooth functions on
$V(M)$. A \textit{Lie bracket} of two vector fields $a,b \in V(M)$ is a vector
field given by the formula:
\[
[a,b]_j(x) = \sum_{k=1}^d \left( a_k \frac{\partial b_j}{\partial x_k} (x) - b_k \frac{\partial a_j}{\partial x_k} (x) \right) .
\]
\end{definition}
\begin{definition}
Let a PDMP be defined by the systems of differential equations
$ x'=g_i(x)$,  $i \in I=\{0,1,...,k\}$, $k \in \mathbb{N}$. We say that the
H{\"o}rmander's condition holds at a point $x$ if vectors
\[
g_2(x)-g_1(x), \dots,g_k(x) - g_1(x), [g_i,g_j](x)_{1 \leq i, j \leq k}, [g_i, [g_j, g_l]](x)_{1 \leq i, j, l \leq k },\dots
\]
 span the space $\mathbb{R}^d$.
\end{definition}

\begin{definition}
Let $n \in \mathbb{N},\ t>0,\ \tau=(\tau_1,\tau_2,\dots,t-\tau_{n-1}-\dots-\tau_{1})$ and
$i=(i_1,\dots,i_n)$ such that for all
$k \in \{1,\dots,\ n-1 \}$ we have  $\tau_k>0$, $i_k \neq i_{k+1} $ and $ i_k \in \{0,1\}$.
A function
\[
\psi_{x,t,i}(\tau):= \pi_{t-\tau_{n-1}-\dots-\tau_{1}}^{i_n} \circ \pi_{\tau_{n-1}}^{i_{n-1}} \circ \dots \circ \pi_{\tau_{1}}^{i_1}
\]
is called a \textit{cumulative flow} along the trajectiories of the flows $\pi^{i_1},\dots,\ \pi^{i_n} $ with starting point $x$.

\end{definition}

\begin{definition}
We say that a point $x \in X$ \textit{communicates} with $y \in X$ if there exist $\ n \in \mathbb{N}$,
$t>0$, $\tau=(\tau_1,\tau_2,\dots,\ t-\tau_{n-1}-\dots-\tau_{1}) $ and $i=(i_1,\dots,i_n)$ such that $\psi_{x,t,i}(\tau)=y$.
\end{definition}
If every two points from the interior of $X$ communicate, we call this property \textit{communication between states of the process}.
If for  $q_0\in X$ there exists $p\in X$ such that $q_0$ communicates with $p$ and the H\"{o}rmander's condition holds at the point $p$, then
$q_0$ satisfies condition (b). This fact is a simple consequence of \cite[Theorem 4]{bakhtin}.

Let us
denote by $a_i(\mathbf {x})$ a vector field representing the system $(\ref{e:2})$ with a fixed value of $i \in \{0,1\}$
at a point $\mathbf {x} \in [0,1]^3.$ After short calculation of the following expressions:
\[
a_1-a_0=(1,0,0),\ [a_0,a_1]=(1,-a,0) ,\ [a_0,[a_0,a_1]]=(1,-(a^2+a),ab);
\]
we obtain three linear independent vectors in $\mathbb{R}^3$ space. Hence, these
vectors span $\mathbb{R}^3$ and condition (b) of Theorem \ref{th:2} holds.
However, it gets more difficult to check condition (a), because it does not hold on the whole space $[0,1]^3 \times \{0,1\}$.
We will prove that (a) holds on $\mathbb{A}$. Moreover, $A$ is a stochastic attractor,
i.e. a measurable subset of $[0,1]^3 $ such that for every density $f \in L^1(\mathbb{A})$ we have
\begin{equation}\label{e:400}
\lim_{t \rightarrow \infty} \int_{\mathbb{A}} P(t) f(\mathbf {x},i)\, dx\, di=\lim_{t \rightarrow \infty} \mathbb{P}(\zeta(t) \in A) = 1.
\end{equation}

First, we show that $\mathbb{A}$ is an invariant set for the process. It follows from the fact that if we take
any $\mathbf {x} \in A,$ then we
stay in $A$ under the action of both semi-flows given by $(\ref{e:2})$ with the initial condition $\mathbf {x}$.
In other words, we check that if we take any $\mathbf {x} \in A$ and $t>0,$ then both $\pi_t^0(\mathbf {x})$ and $\pi_t^1(\mathbf {x})$ stay in $A.$ Since the process switches between these two flows, its trajectories cannot leave $A.$

Let us define the set:

\begin{equation*}
D=\{(x-y+z-w,\ x^a-y^a+z^a-w^a,\ x^b-y^b+z^b-w^b)\colon 1\geqslant x\geqslant y \geqslant z \geqslant w \geqslant 0 \}.
\end{equation*}
We prove that $D$ is the same set as $A.$ For $s \in \{1,a,b\}$ and any
\[
d_0=(x_0-y_0+z_0-w_0,\ x_0^a-y_0^a+z_0^a-w_0^a,\ x_0^b-y_0^b+z_0^b-w_0^b) \in D
\]
we have
\begin{align*}
x_0^s-y_0^s+z_0^s-w_0^s&=x_0^s \left(1-\dfrac{y_0^s}{x_0^s}+\dfrac{z_0^s}{x_0^s}-\dfrac{w_0^s}{x_0^s} \right)
= x_0^s({x'_0}^s-{y'_0}^s+{z'_0}^s)\\
&=(x_0 x'_0)^s-(x_0 y'_0)^s+(x_0 z'_0)^s,
\end{align*}
where the second equality follows from the equivalence of two definitions of the set $A$ (see $\ref{e:500}$).
Therefore, $A=D$. Since it is easy to check that for any $t>0$ and $i \in \{0,1\}$ we have $\pi_t^i(A)\subset D$,
$\mathbb{A}$ is an invariant set for our process.

Using the formula $\ref{e:300},$ there exist $t_1^0>0$, $t_2^0>0$, and $t_3^0>0$ such that we have $\pi_{t_3^0}^0\pi_{t_2^0}^1\pi_{t_1^0}^0 \mathbf 1=(\frac12,\frac12,\frac12)$.
 From continuous dependence of the solutions on the initial condition, we can find $\delta>0$ and $\varepsilon>0$ such that for every
 $\mathbf y\in \mathbb R^3$ with $\|\mathbf y - \mathbf 1\|<\delta$ and $\mathbf t=(t_1,t_2,t_3)$ with $|t_i-t_i^0|<\varepsilon$ for $i=1,2,3$
we have $\pi_{t_3}^0\pi_{t_2}^1\pi_{t_1}^0(\mathbf  y)\in \mathrm{Int\,}A$. Moreover, there exists $t_0>0$ such that
 $\|\pi_t^1(\mathbf x)- \mathbf 1\|<\delta$ for each point $\mathbf x\in [0,1]^3$ and $t>t_0$.
 Thus $\pi_{t_3}^0\pi_{t_2}^1\pi_{t_1}^0\pi_{t_0}^1(\mathbf  x)\in \mathrm{Int\,}A$ if $\mathbf x\in [0,1]^3$
 and $|t_i-t_i^0|<\varepsilon$ for $i=1,2,3$ and $t>t_0$. The probability that the sequence of the consecutive five jump moments $T_n,\dots, T_{n+4}$ has the properties:  $T_{n+1}-T_n>t_0$ and $T_{n+i+1}-T_{n+i}\in (t_i^0-\varepsilon,t_i^0+\varepsilon)$ for $i=1,2,3;$
is bounded from below (see \cite{bobrow}) by some positive number $\eta>0$. Since our PDMP enters the cube $[0,1]^3$ (see the discussion under the formula (\ref{e:sim})) we have also that
it enters the interior of the attractor with probability one, which completes the proof of $(\ref{e:400}).$

To prove the remaining condition (a), we use the fact that it is equivalent to communication
between states for $\mathbf {x},\mathbf {y} \in \operatorname{int}A $
and fixed $i \in \{0,1\}.$ In other words, we show that the cumulative flow
consisting the flows of $\pi^0$ and $\pi^1$ with the initial condition $\mathbf {x} \in \operatorname{int}A $
generates whole $A$.
The question to face is: does the total control of the system $(\ref{e:4})$ between two arbitrary points in the interior of $A$ exist?
The problem of control for linear dynamical systems has been extensively studied in the past (\cite{col}, \cite{kl}),
but this special case appears to be relatively far from classical results of the controllability theory
and seems to not undergo any of those procedures. However, the proof of the communication between states property in this special case is surprisingly simple.
Due to symmetry of $A,$  we consider only these cumulative flows which begin from $\pi^0$. The case when we start from $\pi^1$ is analogous. Fix a point $\mathbf x\in A$. After compositing four transformations we obtain
\[
\psi_{\mathbf x,\mathbf
t,0}:=\pi_{t_4}^1\pi_{t_3}^0\pi_{t_2}^1\pi_{t_1}^0(\mathbf {x})=\mathbf{1} -
\pi_{t_4}^0 \mathbf{1} + \pi_{t_3+t_4}^0 \mathbf{1} - \pi_{t_2+t_3+t_4}^0
\mathbf{1}+ \pi_{t_1+t_2+t_3+t_4}^0(\mathbf {x}) .
\]
Fix $\varepsilon >0$ and take any points $x,y,z$ such that $\varepsilon \le
z\le y\le x\le 1$. Then we can find $t_i\ge 0$, for $i=1,2,3,4$, such that
$x=e^{-t_4}$, $y=e^{-t_3-t_4}$,
$z=e^{-t_2-t_3-t_4}$, and $\varepsilon=e^{-t_1}z$. Thus
\[
\psi_{\mathbf x,\mathbf t,0}
=(1-x+y-z+\varepsilon x_1,1-x^a+y^a-z^a+\varepsilon^a
x_2,1-x^b+y^b-z^b+\varepsilon^b x_3).
\]
Let
\[
A_{\varepsilon}=\{(1-x+y-z,1-x^a+y^a-z^a,1-x^b+y^b-z^b)\colon \, \varepsilon
\le z\le y\le x\le 1\}
\]
and $\mathbf v_{\varepsilon}=(\varepsilon  x_1,\varepsilon^a
x_2,\varepsilon^b x_3)$.
Then, starting from the point $\mathbf x$ and using a composition of four
transformations we communicate with each point from the set
\[
A_{\varepsilon}+\mathbf v_{\varepsilon}:=\{\mathbf y+\mathbf
v_{\varepsilon}\colon \mathbf y\in A_{\varepsilon}\}.
\]
Since $\mathrm{Int\ } A\subset \bigcup\limits_{\varepsilon>0}
A_{\varepsilon},$ we conclude that
we can join $\mathbf {x}$ and any interior point of $A$ by $\psi_{\mathbf
{x},t,0}$ (note that this property may not hold for the boundary points of
$A$).
As a result, condition (a) from Theorem~\ref{th:2} is satisfied and the
semigroup  $\{P(t)\}_{t \geqslant 0} $ is asymptotically stable.

\medskip

{\it Remark 3.} The proof  in the general case (i.e. including the subcases  $a = b$,  $a = 1$ and $b = 1$) is similar to the presented above, but technically more difficult. We do not change the variables in the system (\ref{e:2}) (see Appendix A), but we define the attractor for the process as the closure of the set:
\begin{equation}
A_1=\{\mathbf x= \pi_{s_1}^0 \mathbf{1} -\pi_{s_2}^0 \mathbf{1}
+ \pi_{s_3}^0 \mathbf {1} :\,\, 0\le s_1\le s_2\le s_3\}
\end{equation}
or equivalently as the closure of the set:
\begin{equation}
A_2=\{\mathbf x=\mathbf{1}- \pi_{s_1}^0 \mathbf{1} +\pi_{s_2}^0 \mathbf{1}
- \pi_{s_3}^0\mathbf {1} :\,\, 0\le s_1\le s_2\le s_3\}.
\end{equation}
We have $\textrm{cl\,}A_1=\textrm{cl\,}A_2,$ because both sets $A_1$ and
 $A_2$ have the same boundaries:
 \[
 \textrm{cl\,}\{\mathbf x=\pi_{s_1}^0 \mathbf{1} -\pi_{s_2}^0 \mathbf{1}:\,\, 0\le s_1\le s_2\}\cup
 \textrm{cl\,}\{\mathbf x=\mathbf{1}- \pi_{s_1}^0 \mathbf{1} +\pi_{s_2}^0 \mathbf{1}
 :\,\, 0\le s_1\le s_2\}.
\]
Having these formulas for the attractor one can prove that indeed this set is an attractor and that any two points in the interior of the attractor can communicate.

\bibliography{mybib2}
\bibliographystyle{plain}

\end{document}